\documentclass[a4paper,11pt,oneside,twocolumn]{article}

\usepackage[T1]{fontenc}
\usepackage[utf8]{inputenc}
\usepackage[italian,english]{babel}
\usepackage{abstract}
\usepackage{color}
\usepackage{graphicx}
\usepackage[italian]{varioref}
\usepackage{fancyhdr}
\usepackage{charter}
\usepackage{booktabs}
\usepackage{amsthm}
\usepackage{dsfont}
\usepackage{indentfirst}
\usepackage{amsmath}
\usepackage{amsfonts}
\usepackage{amssymb}
\usepackage{enumitem}
\usepackage{footnote}
\graphicspath{{Immagini/}}

\begin{document}

\title{Corner Detection and Arc Approximation of Planar Data Set}
\author{Maurizio Scarparo}
\date{}
\twocolumn[
\maketitle
\hspace{21pt}\rule{409pt}{.4pt}
\begin{onecolabstract}
The aim of this paper is to present a new method of approximation of planar data set using only arcs or segments. The first problem we are trying to solve is the following: the CNC machines can work only with simple curves (arcs or segments, indeed). So, if we have a contour of an object and we want to cut, for example, glass or metal following that profile, we need a continuous curve composed by linear or circular paths only. Moreover, we want to minimize the number of these paths. The second problem is the following: the contour of an object is detected by a specific laser, which collects a discrete data set. The laser detects a good approximation of a point if this point is not a corner of the object. So in many cases our data set will not contain the real corners. Our task is to present a method on how to find these particular points. A third purpose is to study the regularity of the approximating curve, which can be $G^1$ or $C^0$ continuous. The entire method is developed using single arcs, segments and in certain cases biarcs, in order to ensure the smoothness of the final path.

\vspace{5pt}

\textbf{Contacts}: maurizio.scarparo24@gmail.com  
\end{onecolabstract}
\hspace{21pt}\rule{409pt}{.4pt}\vspace{10pt}
]

\section{Introduction}

The problem of approximation of planar data set has been already studied by several authors, such as Les A. Piegl, Wayne Tiller, Hyungjun Park, D. Meek, D. Walton (see page 17). In their articles these authors explain the biarc model and how to apply it in CNC machining. One of the most powerful methods, introduced by Piegl and Tiller, consists in finding the least square B-spline and approximating this one with biarcs (see [7], [8], [9]). The problem of corner detection has also been studied by authors such as C. Harris, M. Stephens, H. Freeman, L. Davis. However, we decided to analyze these problems in a different way. The first task was to elaborate a new method of approximating a planar data set without using B-splines. Indeed, we found out that the B-spline method is very powerful when the data set has no corner points, but it is quite difficult to get a good approximation when the data set is not smooth. Moreover, we found out that the approximation heavily depends on the data set and on the choice of the tangent vectors, so we wanted to cancel the dependence on the vectors. We also discovered that single arcs are better than biarcs if we want to minimize the number of paths. Secondly, we wanted to develop a straightforward method for corner detection in a planar data set captured by a laser. As already explained in the abstract, this data set may not include the exact corner points, because of the laser errors. This is the main difference between our algorithm and the methods explained in [1], [2], [6], [11], [12]: all methods developed there need the corner points to be included in the data set. We noticed that this problem has been carefully studied for digital images using pixel properties, but as far as we know there are not many articles concerning this topic outside image analysis.


\section{Least square constrained arcs}
\label{sec: first}
The first time we try to approximate a planar data set with only biarcs we can find out an inconvenient: if the distance between two consecutive points is small and the tangent vectors are quite different, the biarc between the two points will be \textit{serpentine}, that means useless and insignificant for a good approximation. So the first step is to forget the biarc model for a while and to focus on the \textit{least square constrained arcs}. Let $\mathbf{P}_0$ \footnote{Notation: $\mathbf{P}_i = [x_i,y_i]$}, $\mathbf{P}_1$, $\dots$, $\mathbf{P}_{N+1}$ be  $N+2$ distinct planar points. We want to find one particular arc with initial point $\mathbf{P}_0$, final point $\mathbf{P}_{N+1}$, approximating $\mathbf{P}_1$, $\dots$, $\mathbf{P}_N$ and minimizing an objective error function. Let $R$, $\mathbf{C}$, be respectively the radius and the center of our arc. The function we want to minimize is the following:

\begin{equation}
\label{eqn: f}
f(t) = \sum_{k = 1}^{N} {\bigg\vert R - \parallel\mathbf{P}_k - \mathbf{C}(t){\parallel} \bigg\vert}^2.
\end{equation}

The variable $t$ represents the distance between the centre $C$ and the segment $\overline{\mathbf{P}_0\mathbf{P}_{N+1}}$. In particular, let $d =\ \parallel\mathbf{P}_{N+1}-\mathbf{P}_0\parallel$ be the distance between $\mathbf{P}_0$ and $\mathbf{P}_{N+1}$, let $\mathbf{u}$ be the normalized vector

\begin{equation}
\mathbf{u} = \frac{\mathbf{P}_{N+1}-\mathbf{P}_0}{d}
\end{equation}

and let ${\mathbf{u}}'$ be the $\frac{\pi}{2}$ counterclockwise rotation of $\mathbf{u}$. Then, the center $\mathbf{C}$ can be represented by the following:

\begin{equation}
\mathbf{C}(t) = \mathbf{P}_0 + \frac{d}{2}\mathbf{u} + t{\mathbf{u}}'.
\end{equation}

Let us discuss about the objective function (1). This is a non negative continous function on $\mathds{R}$. One can think that $f$ depends on two unknown $R$ and $\mathbf{C}$; however, since the initial point of our arc is $\mathbf{P}_0$ and the final one is $\mathbf{P}_{N+1}$, the radius depends on the center and 
\begin{equation}
\label{eqn: radius}
R(t) =\ \parallel \mathbf{P}_0 - \mathbf{C}(t) \parallel\ =\ \parallel \mathbf{P}_{N+1} - \mathbf{C}(t) \parallel.
\end{equation}
So, according to (4), the function $f$ in (1) becomes
\begin{equation}
\begin{split}
f(t) =\ & N{\parallel\mathbf{P}_0 - \mathbf{C}(t) \parallel}^2 + \sum_{k = 1}^N {\parallel \mathbf{P}_k - \mathbf{C}(t) \parallel}^2\\
& - 2\parallel\mathbf{P}_0 - \mathbf{C}(t) \parallel\sum_{k = 1}^N {\parallel \mathbf{P}_k - \mathbf{C}(t) \parallel}.
\end{split}
\end{equation}
According to Weierstrass theorem, $f$ admits minimum on each compact set $[a,b] \subset \mathds{R}$. Moreover, one can easily see that the first derivative of $f$ is always continuous on $\mathds{R}$ if the center $\mathbf{C}(t)$ and the points $\mathbf{P}_k$ are distinct for every $k = 0,1,\dots,N+1$. We will prove that $f$ admits finite limits at $\pm\infty$ and these limits are the same.

\theoremstyle{plain}
\newtheorem{theorem}{Theorem}[section]

\begin{theorem}
\label{thm: one}
$f$ has finite limit at $\pm\infty$ and
\begin{equation}
\label{eqn: mine}
\lim_{t \to \pm\infty} f(t) = \sum_{i=1}^N {\bigg[\mathbf{u}' \cdot {(\mathbf{P}_0 - \mathbf{P}_i)}\bigg]}^2,
\end{equation}
where $\cdot$ means the standard scalar product in ${\mathds{R}}^2$.
\end{theorem}

\proof

We have that
\begin{equation}
\begin{split}
&{\parallel \mathbf{P}_0 - \mathbf{C}(t) \parallel}^2 = {\bigg\| \frac{d}{2}\mathbf{u} + t\mathbf{u}' \bigg\|}^2 = t^2 + \frac{d^2}{4},\\
&{\parallel\mathbf{P}_i-\mathbf{C}(t)\parallel}^2 = t^2+2t\bigg[\mathbf{u}'\cdot\bigg(\mathbf{P}_0-\mathbf{P}_i+\frac{d}{2}\textbf{u}\bigg)\bigg]\\
& + {\bigg\|\mathbf{P}_0-\mathbf{P}_i+\frac{d}{2}\textbf{u}\bigg\|}^2.
\end{split}
\end{equation}
Let $\mathbf{D}_i = \mathbf{P}_0-\mathbf{P}_i+\frac{d}{2}\textbf{u}$. We notice that $f$ is the sum of $N$ functions $f_i(t)$ such that
\[f_i(t) = {(\|\mathbf{P}_0 - \mathbf{C}(t)\| - \|\mathbf{P}_i - \mathbf{C}(t)\|)}^2 = \]
\[ = {(g_0(t) - g_{i}(t))}^2,\]
where
\begin{equation}
\begin{split}
g_0(t) =& \sqrt{t^2 + \frac{d^2}{4}},\\
g_i(t) =& \sqrt{t^2 + 2t\big(\textbf{u}' \cdot {\mathbf{D}_i}\big)+{\|\mathbf{D}_i\|}^2}.
\end{split}
\end{equation}
Multiplying $f_i(t)$ by $\frac{{(g_0(t) + g_i(t))}^2}{{(g_0(t) + g_i(t))}^2}$ one obtains
\begin{equation}
\label{eqn: second}
\begin{split}
& f_i(t) = \frac{{(g_0^2(t) - g_i^2(t))}^2}{{(g_0(t) + g_i(t))}^2} = \\
=& {\frac{{\bigg(\frac{d^2}{4} - 2t\big(\textbf{u}' \cdot {\mathbf{D}_i}\big) - {\|\mathbf{D}_i\|}^2\bigg)}^2}{{\bigg(\sqrt{t^2 + \frac{d^2}{4}} + \sqrt{t^2 + 2t\big(\textbf{u}' \cdot {\mathbf{D}_i}\big)+{\|\mathbf{D}_i\|}^2}\bigg)}^2}}
\end{split}
\end{equation}
By equation (9) one can easily see that
\begin{equation}
\lim_{t \to \pm\infty} f_i(t) = \big(\textbf{u}' \cdot {\mathbf{D}_i}\big)^2 = \big(\textbf{u}' \cdot (\mathbf{P}_0 - \mathbf{P}_i)\big)^2,
\end{equation}
and (6) is proved using limit properties on the sum.
\endproof

\theoremstyle{plain}
\newtheorem{corollary}{Corollary}[section]
\begin{corollary}
The limit of $f$ at $\pm\infty$ is zero if and only if all the points $\mathbf{P}_1,\dots,\mathbf{P}_N$ are collinear on the segment $\overline{\mathbf{P}_0\mathbf{P}_{N+1}}$.
\end{corollary}

Theorem 2.1 is quite important because we can say that, if $L$ is the limit of $f$ at $\pm\infty$,
\[ \forall\ \varepsilon > 0\ \exists\ \overline{t} \in \mathds{R} : |t| > |\overline{t}|\Longrightarrow |f(t) - L| < \varepsilon. \]
So the set $B(0,|\overline{t}|) = \{t \in \mathds{R} : |t| \leq |\overline{t}|\}$ is compact, $f$ is continuous on $B(0,|\overline{t}|)$ and admits minimum there. 
In order to find a minimum of $f$ one can apply one of the optimization techniques, such as Newton or Quasi-Newton methods, BFGS methods and others. In literature, one can find many methods of building a least square circle given some planar points. One of the most efficient methods is Taubin's method, explained in detail in [15]. The method explained int this paper has one main difference with Taubin's: the least square arc is constrained, that is it passes through two given points. Taubin's arc is better than ours, because it minimizes the error on all the possible centers and radii, but in CNC approximation it's very difficult to find two consecutive Taubin's arcs connecting each other in a continuous way.

In the following page, there is a simple example of the function $f$ and the related least square arc. 

\section{Searching longest arcs}

In this section we want to explain how to apply least square constrained arcs in the approximation of a planar data set without corners (we suppose that the profile of our object is locally smooth). Suppose that $\mathbf{P}_1,\mathbf{P}_2,\dots,\mathbf{P}_N$ are $N$ consecutive points approximating one section of the object contour. 

\begin{figure*}[htbp]
\centering
\includegraphics[scale=0.55]{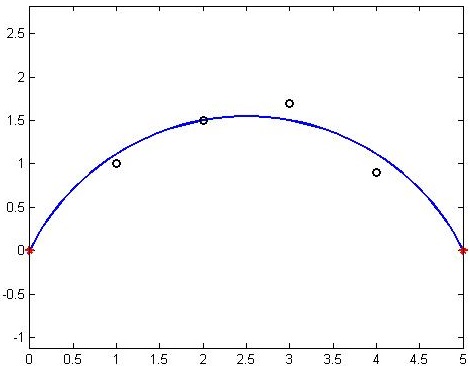}%
\qquad\qquad
\includegraphics[scale=0.55]{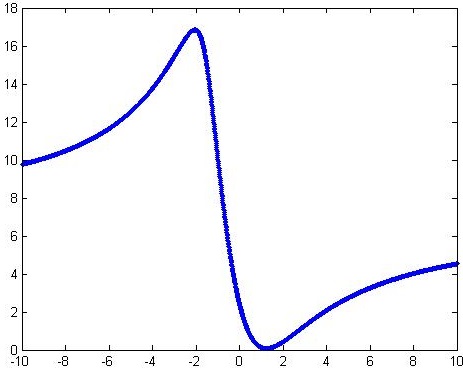}
\caption{Example of least square constrained arc passing through (0,0),\ (5,0) (on the left) and the relative function $f$ (on the right).}
\end{figure*}
The first step of the method consists in trying to approximate the first four points with a constrained arc (it starts from four because we already know that there always exists one arc passing through three points, eventually with infinite radius). We decide an error limit and we compute the distance of each point $\mathbf{P}_k$ from the arc. If there exists one point whose distance is greater than our limit, we decrease the number of approximated points and, in this case, we build the arc passing through the first three points. Conversely, if all the points are within the error limit, we increase the number of approximated points and we try to build a longer arc. The method repeats these steps for each point $\mathbf{P}_k, k = 1,2,\dots,N-2$. At the end we obtain a $N \times 1$ matrix $M$: in the $k^{th}$ place there will be the number of points that can be approximated with a single least square arc starting from $\mathbf{P}_k$ --- included the first and the last points ---. We will set $M(N-1) = M(N) = 0$. Here follows the algorithm (which needs at least four points in order to be applied):

\hrulefill
\begin{enumerate}
\item Set $i = 1,\ k = 3,\ flag = 1$ and choose $\varepsilon > 0$;
\item \textbf{While} $i \leq N-2$ \textbf{do}
\addtolength{\itemindent}{1cm}
\item Set $l = k - 1$;
\item Find the center $\mathbf{C}$ and the radius $R$

\hspace{25pt} of the least square arc passing

\hspace{25pt} through $\mathbf{P}_i$, $\mathbf{P}_{i+k}$;
\item \textbf{For} $j = 1,\dots,l$
\addtolength{\itemindent}{1cm}
\item $d =\ \mid R - \parallel \mathbf{P}_{i+j} - \mathbf{C} \parallel \mid$;
\item \textbf{If} $d > \varepsilon$ set $flag = 0$; \textbf{break};
\addtolength{\itemindent}{-1cm}
\item \textbf{If} $flag = 0$ 
\addtolength{\itemindent}{1cm}
\item build the least square arc 

\hspace{53pt} passing through $\mathbf{P}_i$, $\mathbf{P}_{i+k-1}$
\item $M(i) = k$, $flag = 1$, $i = i + 1$;
\item \textbf{If} $i = N-2$ set $k = 2$

\hspace{53pt} \textbf{Else} set $k = 3$;
\addtolength{\itemindent}{-1cm}

\hspace{25pt} \textbf{Else}
\addtolength{\itemindent}{1cm}
\item $k = k + 1$;
\item \textbf{If} $i + k > N$ 
\addtolength{\itemindent}{1cm}
\item build the least square

\hspace{82pt} arc passing through 

\hspace{82pt} $\mathbf{P}_i$, $\mathbf{P}_{i+k-1}$
\item $M(i) = k$, $i = i + 1$;
\item \textbf{If} $i = N-2$ set $k = 2$

\hspace{82pt} \textbf{Else} set $k = 3$;
\addtolength{\itemindent}{-3cm}
\item $M(N-1) = 0,\ M(N) = 0$.
\end{enumerate}

\hrulefill

One can notice that $M(k) \geq 3$ for all $k = 1,\dots,N-2$, so that the worst case occurs when all the elements of $M$ are equal to three. Next step is to decide how to choose the longest arcs.

\section{Selection of longest arcs}

We will explain the choice of longest arcs with two simple examples. Fix $k \geq 1$ and suppose that $M(k) = n$, for a certain $n \geq 3 \in \mathds{N}$. Suppose also that $M(i) < M(k)$ for all $i = k+1,k+2,\dots,k+n-1$. This is the simplest case: the arc starting from $\mathbf{P}_k$ is approximating an higher number of points than the $n-1$ consecutive arcs, so it is good for an optimized approximation. Now suppose that an index $j$ exists such that $ 1 \leq j \leq n-1$ and $M(k+j) > M(k)$. In this situation we can't say that the arc starting from $\mathbf{P}_k$ is the best for our task. Let 
\[S = \bigg\lfloor \frac{M(k)}{2} \bigg\rfloor\]
and let $I_k$ be the set
\[I_k = \{j \in \mathds{N}\ \vert\ k+1 \leq j \leq k + S\}.\]
Then, we decide to apply the following test: if $M(j) \leq M(k)$ for all $j \in I_k$, the arc starting from $\mathbf{P}_k$ is good for our approximation and we pass from $\mathbf{P}_k$ to $\mathbf{P}_{k+M(k)-1}$ with a single arc. We repeat the process on $\mathbf{P}_{k+M(k)-1}$ (the new starting point). Conversely, if an index $j \in I_k$ exists such that $M(j) > M(k)$ we pass directly from the point $\mathbf{P}_k$ to $\mathbf{P}_j$ without any segments or arcs --- for the moment --- and we apply again the test on the new starting point $\mathbf{P}_j$. This kind of selection is justified by the following fact: suppose there exists an index $j \notin I_k$ such that $M(j) > M(k)$ and $M(k)$ is a high number ($M(k) \geq 20$ for example ); then, the gap between $\mathbf{P}_k$ and $\mathbf{P}_j$ can be so wide that we are not able to obtain a good approximation for all those points among them, so the test is applied only on the set $I_k$.

Here follows the algorithm ($N$ is, as usual, the number of points we are approximating):

\hrulefill
\begin{enumerate}

\item Set $i = 1$, $G = M$;
\item \textbf{While} $i < N$ \textbf{do}
\addtolength{\itemindent}{1cm}
\item $S = \bigg\lfloor \frac{M(i)}{2} \bigg\rfloor$, $j = i+1$, $flag = 1$,

\hspace{25pt} $e = i + M(i) - 1$;
\item \textbf{While} $j \leq i + S$ \textbf{and} $flag = 1$ \textbf{do}
\addtolength{\itemindent}{1cm}
\item \textbf{If} $M(j) \leq M(i)$ set $j = j + 1$

\hspace{53pt} \textbf{Else} 
\addtolength{\itemindent}{1cm}
\hspace{53pt} \item Set $flag = 0$;
\hspace{53pt} \item \textbf{for} $k = i,\dots,j-1$ 

\hspace{82pt} set $G(k) = 0$;
\hspace{53pt} \item Set $i = j$;

\hspace{82pt} $e = i + M(i) - 1$; \textbf{break};
\addtolength{\itemindent}{-2cm}
\item \textbf{If} $flag = 1$
\addtolength{\itemindent}{1cm}
\item \textbf{For} $k = i+1,\dots,e-1$

\hspace{53pt} $G(k) = 0$;
\item Set $i = e$;

\end{enumerate}
\hrulefill

In this algorithm we introduce a new matrix $G$. At the beginning, $G$ is equal to $M$, but during the steps the element $G(i)$ may become null. This happens in two cases, if $\mathbf{P}_k$ is a starting point:

\begin{itemize}
\item $M(k) \geq M(j)$ for all $j \in I_k$. Then we set $G(j) = 0$ for all $j$ such that $k+1 \leq j \leq k + M(k) - 2$;
\item there exists some index $j \in I_k$ such that $M(k) < M(j)$. Then we set $G(i) = 0$ for all $i \in I_k$ satisfying $k + 1 \leq i \leq j-1$ and for $i = k$.
\end{itemize}

At the end of the process we obtain a $N\times 1$ matrix $G$ with elements $G(k)$ satisfying:
\renewcommand{\labelitemi}{$\diamond$}
\begin{itemize}
\item $G(k) > 0$, if there exists a single arc passing through the points $\mathbf{P}_k$, $\mathbf{P}_{k+G(k)-1}$ and approximating all the other points among them, i.e. if the point $\mathbf{P}_k$ is a starting point for one arc;
\item $G(k) = 0$, if the point $\mathbf{P}_k$ is not a starting point, so $G(N-1) = G(N) = 0$.
\end{itemize}
\renewcommand{\labelitemi}{$\bullet$}
Suppose now $G(k) > 0$ and $G(k+G(k)-1) > 0$: then there exist two consecutive arcs (continous arcs) with common point $\mathbf{P}_{k+G(k)-1}$. This is the best case: all points between $\mathbf{P}_k$ and $\mathbf{P}_{k+G(k)-1+G(k+G(k)-1)-1}$ are approximated within the tolerance. But what happens when two consecutive arcs don't have a common point? This is the next step: how to fill the gaps between the built arcs.

\section{Constrained Biarcs}
\theoremstyle{definition}
\newtheorem{definition}{Definition}[section]
Before proceeding to the next step, we want to discuss about biarcs and their applications.
\begin{definition}
\label{def: biarc}
We say that two circular arcs $a_1$, $a_2$ form a constrained biarc passing through given points $\mathbf{P}_s$, $\mathbf{P}_e$ with unit tangent vectors $\mathbf{T}_s$, $\mathbf{T}_e$, if they have a common point $\mathbf{P}_2$ satisfying the following properties:
\begin{itemize}
\renewcommand{\labelitemi}{$\ast$}
\item $a_1$ starts from $\mathbf{P}_s$, ends in $\mathbf{P}_2$ and $\mathbf{T}_s$ is tangent to $a_1$ in $\mathbf{P}_1$, with orientation corresponding to a parametrization of $a_1$ from $\mathbf{P}_s$ to $\mathbf{P}_2$;
\item $a_2$ starts from $\mathbf{P}_2$, ends in $\mathbf{P}_e$ and $\mathbf{T}_e$ is tangent to $a_2$ in $\mathbf{P}_e$, with orientation corresponding to a parametrization of $a_2$ from $\mathbf{P}_2$ to $\mathbf{P}_e$;
\item $a_1$ and $a_2$ have proportional tangent vectors in $\mathbf{P}_2$ (i.e. the biarc is $G^1$ continous in $\mathbf{P}_2$).
\renewcommand{\labelitemi}{$\bullet$}
\end{itemize}
\end{definition}
In this paper we will use the standard notation used by Piegl and Tiller in [7], [8]. The control points of the biarc are determined by
\[\mathbf{P}_1 = \mathbf{P}_s + \alpha\mathbf{T}_s,\]
\[\mathbf{P}_3 = \mathbf{P}_e - \beta\mathbf{T}_e,\]
with both $\alpha$ and $\beta$ positive. Using this convention, the sweep angle of each arc of the biarc will always be in $[0,\pi)$.
Now, let $\mathbf{V} = \mathbf{P}_e - \mathbf{P}_s$. We will prove that there always exists a biarc passing through $\mathbf{P}_s$, $\mathbf{P}_e$ when
\begin{equation}
\label{eqn: cond1}
\mathbf{T}_s\cdot\mathbf{T}_e \neq 1
\end{equation}
and
\begin{equation}
\label{eqn: cond2}
\mathbf{V}\cdot\mathbf{V} \neq \frac{2(\mathbf{V}\cdot\mathbf{T}_s)(\mathbf{V}\cdot\mathbf{T}_e)}{\mathbf{T}_s\cdot\mathbf{T}_e-1}.
\end{equation}

\begin{theorem}
Given two distinct planar points $\mathbf{P}_s$, $\mathbf{P}_e$ with associated unit vectors $\mathbf{T}_s$, $\mathbf{T}_e$ satisfying (11) and (12), there always exists a biarc passing through $\mathbf{P}_s$, $\mathbf{P}_e$ with tangents $\mathbf{T}_s$, $\mathbf{T}_e$.
\end{theorem}
\proof
Let $\mathbf{P}_1 = \mathbf{P}_s + \alpha\mathbf{T}_s$, with $\alpha > 0$. Then it's sufficient to prove that there exists a point $\mathbf{P}_3$ such that $\mathbf{P}_3 = \mathbf{P}_e - \beta\mathbf{T}_e$, $\beta > 0$, and $\parallel \mathbf{P}_3 - \mathbf{P}_1 \parallel = \alpha + \beta$. This is equivalent to the condition \[(\mathbf{P}_3 - \mathbf{P}_1)\cdot(\mathbf{P}_3 - \mathbf{P}_1) = {(\alpha + \beta)}^2.\] Substituting the values of $\mathbf{P}_1$, $\mathbf{P}_3$, we find the condition \[(\mathbf{V}-\alpha\mathbf{T}_s-\beta\mathbf{T}_e)\cdot(\mathbf{V}-\alpha\mathbf{T}_s-\beta\mathbf{T}_e) = {(\alpha + \beta)}^2.\] Expliciting the products, we find \[\mathbf{V}\cdot\mathbf{V}-2\alpha\mathbf{V}\cdot\mathbf{T}_s-2\beta\mathbf{V}\cdot\mathbf{T}_e+2\alpha\beta(\mathbf{T}_s\cdot\mathbf{T}_e - 1) = 0.\] So, expliciting $\beta$:
\begin{equation}
\label{eqn: beta}
\beta = \frac{2\alpha\mathbf{V}\cdot\mathbf{T}_s - \mathbf{V}\cdot\mathbf{V}}{2\alpha(\mathbf{T}_s\cdot\mathbf{T}_e-1) - 2\mathbf{V}\cdot\mathbf{T}_e}.
\end{equation}
Now, assume that exists $\alpha > 0$ such that $\beta > 0$. Then, one can define the junction point $\mathbf{P}_2$ as
\[\mathbf{P}_2 = \frac{\beta}{\alpha+\beta}\mathbf{P}_1+\frac{\alpha}{\alpha+\beta}\mathbf{P}_3.\] Since $\parallel \mathbf{P}_2 - \mathbf{P}_1 \parallel = \alpha$ and $\parallel \mathbf{P}_3 - \mathbf{P}_2 \parallel = \beta$, we can build two arcs $c_1$, $c_2$ satisfying the conditions of the Definition 5.1.
\endproof
The proof of the theorem uses the following
\theoremstyle{plain}
\newtheorem{lemma}{Lemma}[section]
\begin{lemma}
Given two distinct planar points $\mathbf{P}_s$, $\mathbf{P}_e$ with associated unit vectors $\mathbf{T}_s$, $\mathbf{T}_e$ satisfying [11] and [12], there always exists a positive value of $\alpha$ such that $\beta > 0$, where $\beta$ is defined in equation [13].
\end{lemma}
\proof
Let
\[\alpha_1 = \frac{\mathbf{V}\cdot\mathbf{V}}{2\mathbf{V}\cdot\mathbf{T}_s},\ \alpha_2 = \frac{\mathbf{V}\cdot\mathbf{T}_e}{\mathbf{T}_s\cdot\mathbf{T}_e-1},\] \[K = \frac{\mathbf{V}\cdot\mathbf{T}_s}{\mathbf{T}_s\cdot\mathbf{T}_e-1}.\] Then, by hypothesis, $\alpha_1 \neq \alpha_2$. From equation (13) we can write
\[\beta = K\frac{\alpha - \alpha_1}{\alpha - \alpha_2}.\] There are three cases to be examined.
\begin{itemize}
\item \underline{K > 0}. We have $\mathbf{V}\cdot\mathbf{T}_s < 0$, so $\alpha_1 < 0$ and $\beta$ is positive for $\alpha > \text{max}\{\alpha_2,0\}$. Notice that if $\alpha_1 = \alpha_2$ one obtains $\beta \equiv K > 0$ regardless $\alpha$. In this case one arc of the biarc is degenerate (a segment), but the theorem doesn't concern this situation.
\item \underline{K < 0}. We have $\mathbf{V}\cdot\mathbf{T}_s > 0$, so $\alpha_1 > 0$. There are three subcases:
\begin{itemize}
\item $\alpha_2 \leq 0 < \alpha_1$. We choose $\alpha \in (0,\alpha_1)$;
\item $0 < \alpha_1 < \alpha_2$. We choose $\alpha \in (\alpha_1,\alpha_2)$;
\item $0 \leq \alpha_2 < \alpha_1$. We choose $\alpha \in (\alpha_2,\alpha_1)$.
\end{itemize}
\item \underline{K = 0}. We have $\mathbf{V}\cdot\mathbf{T}_s = 0$, so
\[\beta = \frac{-\mathbf{V}\cdot\mathbf{V}}{2(\mathbf{T}_s\cdot\mathbf{T}_e-1)(\alpha-\alpha_2)}.\] We choose $\alpha > \text{max}\{\alpha_2,0\}$.
\end{itemize}
\endproof
The condition of the lemma 
\begin{equation}
\label{eqn: condition}
\mathbf{V}\cdot\mathbf{V} \neq \frac{2(\mathbf{V}\cdot\mathbf{T}_s)(\mathbf{V}\cdot\mathbf{T}_e)}{\mathbf{T}_s\cdot\mathbf{T}_e-1}
\end{equation}
can be read in a geometrical way. Let $\theta_1$, $\theta_2$, $\theta_3$ be respectively the angles between $(\mathbf{V},\mathbf{T}_s)$, $(\mathbf{V},\mathbf{T}_e)$, $(\mathbf{T}_s,\mathbf{T}_e)$. Then, from (14), we find
\[{\parallel \mathbf{V} \parallel}^2 \neq 2\frac{{\parallel \mathbf{V} \parallel}^2\cos{\theta_1}\cos{\theta_2}}{\cos{\theta_3}-1},\]
or, simplifying,
\begin{equation}
\label{eqn: condition1}
\cos{\theta_3} \neq 2\cos{\theta_1}\cos{\theta_2} + 1.
\end{equation}
It's easy to see that a necessary condition to have $\cos{\theta_3} = 2\cos{\theta_1}\cos{\theta_2} + 1$ is that $-1 \leq \cos{\theta_1}\cos{\theta_2} \leq 0$. The condition (15) includes many cases, but in all of these situations one can solve the connection problem using two consecutive biarcs. 
Clearly, things get simpler when $\alpha = \beta$. In this situation, there are only two cases in which we can't build a single arc: $\mathbf{T}_s\cdot\mathbf{T}_e = 1$ and $\mathbf{V}\cdot(\mathbf{T}_s+\mathbf{T}_e) = 0$. In both cases we have to build two consecutive biarcs, i.e. four arcs, matching with $G^1$ continuity. In general, one can define the biarc ratio $r = \frac{\alpha}{\beta}$ instead of choosing a value of $\alpha$ to have $\beta > 0$.

\begin{figure*}[htbp]
\centering
\includegraphics[scale=0.5]{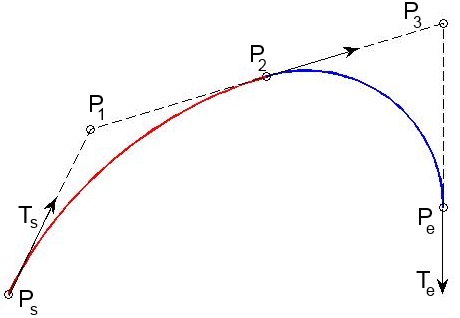}%
\qquad
\includegraphics[scale=0.5]{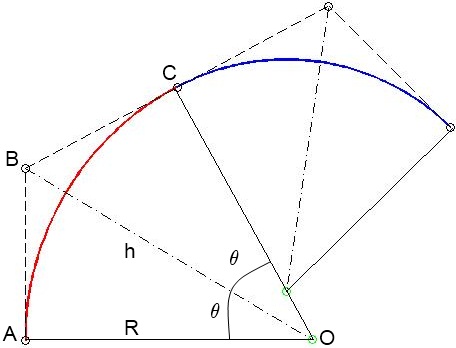}%
\caption{Example of constrained biarc passing through $\mathbf{P}_s = (0,0)$, $\mathbf{P}_e = (5,1)$ with tangent vectors $\mathbf{T}_s = (1,1)$, $\mathbf{T}_e = (0,-1)$ (on the left) and construction of a biarc (on the right).}
\label{fig: esempio}
\end{figure*}

The ratio $r$ is very useful when we want to build a constrained biarc minimizing the distance from given planar points. An efficient algorithm has been studied by Hyungjun Park in [5]: it's called \textit{optimal single biarc fitting}, and its approach is to search an optimized value of $r$ decreasing the width of the range of $r$ in each step. We will use this approach in the next steps.

\section{Biarcs in Bézier form}

Let us consider the biarc in Figure 2 (on the right). Let $\mathbf{ab} = \frac{\mathbf{B}-\mathbf{A}}{\|\mathbf{B}-\mathbf{A}\|}$, $\mathbf{bc} = \frac{\mathbf{C}-\mathbf{B}}{\|\mathbf{C}-\mathbf{B}\|}$. Then we have
\[D = (\mathbf{ab}, \mathbf{bc}) = \cos{2\theta},\]
where $\theta$ is the \textit{half sweep angle} of the biarc, and we suppose that $\theta \in \big[0,\frac{\pi}{2}\big)$. By bisection formulas one can find
\begin{equation}
\begin{split}
\cos{\theta} = \sqrt{\frac{1 + D}{2}},& \ \sin{\theta} = \sqrt{\frac{1 - D}{2}},\\
\tan{\theta} =& \sqrt{\frac{1 - D}{1 + D}}.
\end{split}
\end{equation}
So the explicit formula for $\theta$ becomes
\[\theta = \arctan{\sqrt{\frac{1 - D}{1 + D}}}.\]
Since $R = h\cos{\theta}$ and  $\|\mathbf{B}-\mathbf{A}\| = h\sin{\theta}$ we can find the radius and the centre of the first arc:
\begin{equation}
\begin{split}
R =& \sqrt{\frac{1 + D}{1 - D}}\|\mathbf{B}-\mathbf{A}\|,\\
\mathbf{O} =& \mathbf{B} + \frac{R}{\cos{\theta}}\frac{\mathbf{ab} - \mathbf{bc}}{\|\mathbf{ab} - \mathbf{bc}\|}.
\end{split}
\end{equation}
Suppose $\mathbf{P}_0 = \mathbf{A}$, $\mathbf{P}_1 = \mathbf{B}$, $\mathbf{P}_2 = \mathbf{C}$ and $\{w_0,w_1,w_2\} = \{1, \cos{\theta}, 1\}$. Then we easily obtain the parametrization of the first arc:
\begin{equation}
\mathbf{G}(t) = \frac{\sum\limits_{i = 0}^2 w_i \mathbf{P}_i B_i^2(t)}{\sum\limits_{i=0}^2 w_i B_i^2(t)},\ t \in [0,1],
\end{equation}
where $B_i^n(t)$, $i = 0,\dots,n$, are the $(n+1)$ Bernstein polynomials defined by
\[ B_i^n(t) = {n \choose i} t^i {(1-t)}^{(n-i)},\ t \in [0,1].\]
Now remember that a biarc has five control points $\mathbf{P}_0,\dots,\mathbf{P}_4$, where the first and the last ones are the initial and final point of the biarc, $\mathbf{P}_2$ is the junction point. Let us consider the characteristic functions $\chi_{[0,1]}$, $\chi_{(1,2]}$. Then we can express the biarc in Bézier form:
\begin{equation}
\begin{split}
&\mathbf{B}(t) =\ \chi_{[0,1]}(2t)\frac{\sum\limits_{i = 0}^2 w_i \mathbf{P}_i B_i^2(2t)}{\sum\limits_{i=0}^2 w_i B_i^2(2t)} +\\
&+ \chi_{(1,2]}(2t)\frac{\sum\limits_{i = 2}^4 w_i \mathbf{P}_i B_{i-2}^2(2t-1)}{\sum\limits_{i=2}^4 w_i B_{i-2}^2(2t-1)},\ t \in [0,1],
\end{split}
\end{equation}
where the weight vector is
\[ \{w_i\}_{i=0}^4 = \{1,\cos{\theta_1}, 1, \cos{\theta_2}, 1\}\]
and $\theta_1$, $\theta_2$ are the half sweep angles of the two arcs.

In general, we can write a PCC continous curve as a finite sum of arcs using Bézier form. Let $T_i f(t) = f(t-i)$ be the translation of $f$. Then a PCC curve $C$ composed by $n+1$ arcs has the form
\begin{equation}
C(t) = \left\{\begin{array}{l l} \sum\limits_{i = 0}^n T_i\chi_{[0,1)}(t)T_i G_i(t), & \text{if $t \in [0,n)$}\\
G_n(1), & \text{if $t = n$}
\end{array},\right.
\end{equation} 

where
\begin{equation}
G_i(t) = \frac{\sum\limits_{j = 0}^2 w_{ij}\mathbf{P}_{ij}B_j^2(t)}{\sum\limits_{j = 0}^2 w_{ij}B_j^2(t)}, t \in [0,1],
\end{equation}

is a classical Bézier rational curve, $w_{ij}$ are determined in order to obtain arcs (cosine of the half sweep angles) and $\mathbf{P}_{ij}$ are control points of the arc $i$ (i.e. $G_i(t)$ is an arc with starting point $\mathbf{P}_{i0}$ and ending point $\mathbf{P}_{i2}$). Moreover, to achieve continuity, the control points must satisfy $\mathbf{P}_{i,0} \equiv \mathbf{P}_{i-1,2}$, $i = 1,2,\dots,n$. 

\section{Building the longest arcs}

Recall that in the last step we found the matrix $G$: $G(k) > 0$ if $\mathbf{P}_k$ is a starting point for some arc, $G(k) = 0$ otherwise. Now we want to build the longest arcs using $G$ and to understand the regularity of our final curve. This is quite simple to do when $G(k) > 0$ and $G(k + G(k) - 1) > 0$. As we explained before, this is the case of two consecutive arcs matching with $C^0$ continuity. So we build the arc starting from $\mathbf{P}_k$, ending in $\mathbf{P}_{k+G(k)-1}$ and minimizing the error function $f$ mentioned in section 2 on page 2. Then, we can progress to the next point with the same method, because we already know that the arc starting from $\mathbf{P}_{k+G(k)-1}$ exists. Generally, it's quite difficult that two consecutive longest arcs connect each others with $C^0$ continuity. So, the idea is to proceed in a  different way according to the number of not approximated points between two consecutive arcs. We will explain better what we mean with two examples. Suppose that one of the longest arcs has its ending point in $\mathbf{P}_k$ and that the consecutive arc starts from $\mathbf{P}_{k+1}$. In this case we can build a segment between the two points $\mathbf{P}_k$, $\mathbf{P}_{k+1}$ in order to have $C^0$ continuity. Suppose now that an ending point is $\mathbf{P}_{k}$ but the next starting point is $\mathbf{P}_{k+2}$. In this case, to achieve continuity, we build the arc passing through $\mathbf{P}_k$, $\mathbf{P}_{k+1}$ and $\mathbf{P}_{k+2}$, which always exists. We find the worst case when the number of not approximated points is greater than $3$. We can solve this problem in two different ways: 
\begin{itemize}
\item using constrained biarcs and optimal single biarc fitting;
\item reapplying the algorithms to this smaller section in order to find a new $G$ and new approximating arcs. 
\end{itemize}
Generally, the first approach is not ideal for an optimized curve. In fact, as we underlined at the beginning, the biarc approximation fails sometimes when the section to be approximated is quite small. However, using biarcs ensures $G^1$ continuity, so the first method creates more regular curves. Here follows the algorithm. The matrix $G$ is supposed to be $N \times 1$, with the last two elements equal to zero.

\hrulefill
\begin{enumerate}
\item Set $i = 1$, $k = 0$;
\item \textbf{While} $i \leq N$ \textbf{do}
\addtolength{\itemindent}{1cm}
\item \textbf{If} $G(i) = 0$ \textbf{do} 4, 5
\addtolength{\itemindent}{1cm}
\item $i = i + 1$; $k = k + 1$;
\item \textbf{If} $ i > N$
\addtolength{\itemindent}{1cm}
\item $i = i - 1$; $k = k - 1$;
\item \textbf{If} $k = 1$

\hspace{76pt} build the segment

\hspace{76pt} $\overline{\mathbf{P}_{i-1}\mathbf{P}_i}$;

\hspace{76pt} $k = 0$;

\item \textbf{If} $k = 2$

\hspace{76pt} build the circle passing 

\hspace{77pt} through $\mathbf{P}_{i-2}$, $\mathbf{P}_{i-1}$, $\mathbf{P}_{i}$

\hspace{77pt} $k = 0$;

\item \textbf{If} $ k \geq 3$

\hspace{76pt} reapply the algorithms 

\hspace{75pt} on the set 

\hspace{75pt} $\{\mathbf{P}_{i-k+1},\dots,\mathbf{P}_i\}$ 

\hspace{75pt} or try to approximate 

\hspace{76pt} with biarcs;

\hspace{76pt} $k = 0$;

\item $i = i+1$;

\addtolength{\itemindent}{-2cm}

\item \textbf{Else} \textbf{do} 12, 13, 14, 15

\addtolength{\itemindent}{1cm}

\item \textbf{If} $k = 1$

\hspace{54pt} build the segment $\overline{\mathbf{P}_{i-1}\mathbf{P}_i}$;

\hspace{54pt} $k = 0$;

\item \textbf{If} $k = 2$

\hspace{54pt} build the circle passing 

\hspace{55pt} through $\mathbf{P}_{i-2}$, $\mathbf{P}_{i-1}$, $\mathbf{P}_{i}$

\hspace{55pt} $k = 0$;

\item \textbf{If} $ k \geq 3$

\hspace{53pt} reapply the algorithms on 

\hspace{53pt} the set $\{\mathbf{P}_{i-k+1},\dots,\mathbf{P}_i\}$ 

\hspace{53pt} or try to approximate 

\hspace{54pt} with biarcs;

\hspace{54pt} $k = 0$;

\item $i = i + G(i) - 1$;
\end{enumerate}
\hrulefill

\begin{figure*}[htbp]
\centering
\includegraphics[scale=0.4]{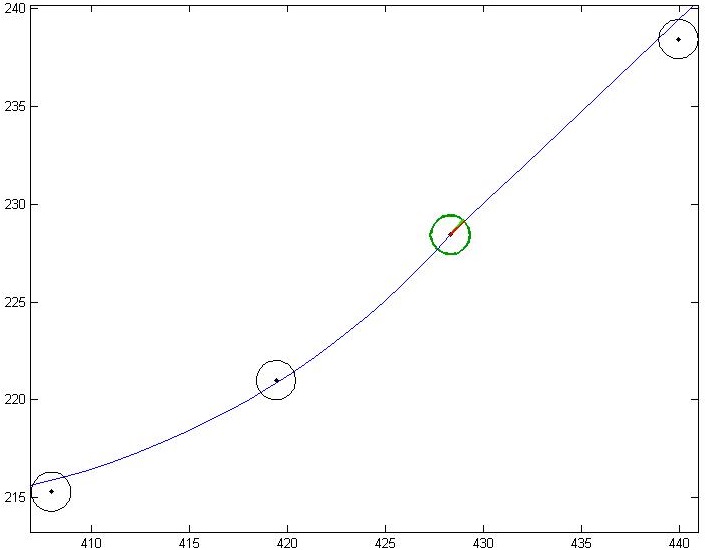}%
\qquad
\includegraphics[scale=0.4]{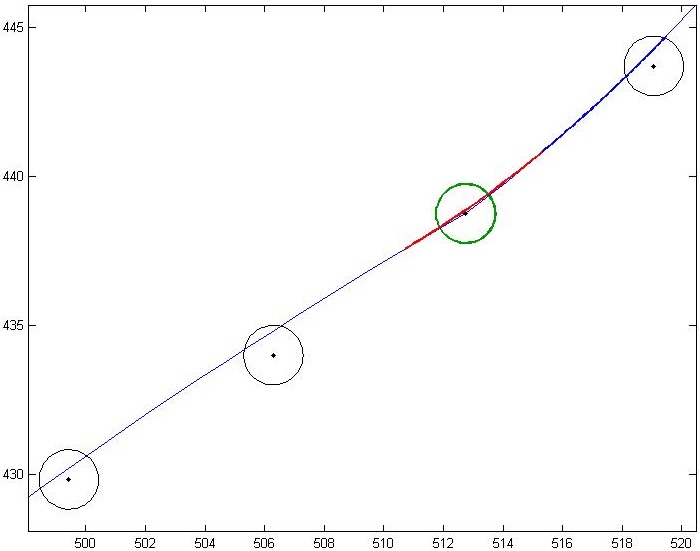}%
\caption{Example of good (on the left) and bad (on the right) junction points (green points). On the left the curve will have $C^0$ continuity. On the right, the red and blue curve is an example of smoothing biarc. Notice that the biarc is within the tolerance and that the final curve will be $G^1$ continuous in the junction point.}
\label{fig: esempi}
\end{figure*}

\begin{figure*}[htbp]
\centering
\includegraphics[scale=0.42]{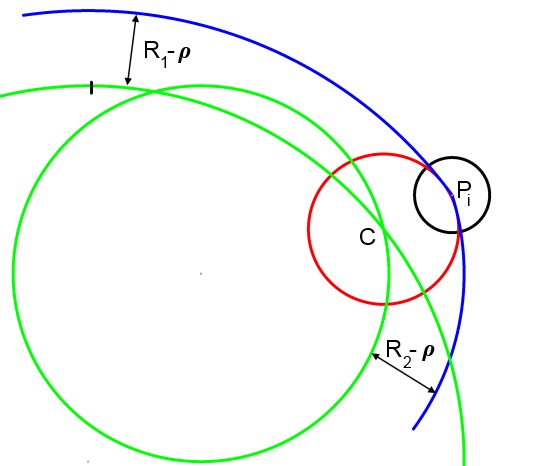}%
\caption{Example of smoothing arc in a neighbourhood of a bad point.}
\label{fig: spiana}
\end{figure*}

\section{Regularity of the approximation}
\label{sec: regular}
Before studying the regularity of the final curve, we need a brief summary. We start from a set of planar data points and we try to build all the possible arcs approximating any number of points (construction of matrix $M$). Then, we make a selection: we take in consideration only the longest arcs (construction of matrix $G$). Finally, if there are some points which haven't been approximated yet, we reapply all the previous steps in the missing sections, creating a continuous curve composed only by arcs or segments.

Clearly, it's quite difficult that all the section of the final curve join each other with $C^1$ continuity, because the approximation uses both arcs and segments. So the next problem we want to discuss is about how to smooth our junction points (j.p.). To do this we first need to recognize \textit{good} j.p. from \textit{bad} j.p. Let us fix a positive tolerance $\varepsilon$ and let $\mathbf{P}_k$ be a junction point for some $k$. Suppose that $\mathbf{P}_k$ is the j.p. of two consecutive arcs $A_{k-1}$, $A_k$ and let $\mathbf{T}_{e,k-1}$, $\mathbf{T}_{s,k}$ be respectively the final unit tangent vector of $A_{k-1}$ and the initial unit tangent vector of $A_k$. We say that $\mathbf{P}_k$ is a \textit{bad} j.p. if

\begin{equation}
\label{eqn: bad}
\mathbf{T}_{s,k} \cdot \mathbf{T}_{e,k-1} = \cos{\vartheta} < \varepsilon,
\end{equation}

where $\vartheta \in [0,\pi]$ is the angle between the two tangent vectors. For example, if $\varepsilon = 0.995$, each angle greater than 6 degrees will be a \textit{bad} angle. Each angle that doesn't satisfy (22) is called a \textit{good} angle. Our task is to smooth out \textit{bad} angles: we can only accept piecewise circular curves (PCC) whose junction points are \textit{good}. This is a practical motivation: if the final object, whose profile is defined by our curve, has \textit{good} junction points, it's easier to smooth out the entire contour, because the angles in good j.p. will be small enough. So if $\mathbf{P}_k$ is a \textit{good} j.p., we don't apply any changes to the PCC (see Figure 3 on the left). Otherwise, in order to increase the regularity of the curve in that particular point, we use biarcs. Indeed, as we have already explained, biarcs are $G^1$ continuous, so they can regularize quite easily a $C^0$ continuous curve. 

Suppose $\mathbf{P}_k$ is a j.p. of two consecutive arcs $A_{k-1}$, $A_k$ and let $\delta$ be a positive fixed tolerance. Then we can find two points $\mathbf{Q}_{k-1}$, $\mathbf{Q}_{k}$ respectively on the arcs $A_{k-1}$, $A_k$ satisfying $\|\mathbf{Q}_{k-1} - \mathbf{P}_k\| = \|\mathbf{Q}_{k} - \mathbf{P}_{k}\| = \delta$, i.e. $\mathbf{Q}_{k-1}$, $\mathbf{Q}_{k}$ are on the circumference with center $\mathbf{P}_{k}$ and radius $\delta$. The idea is to employ these two points as initial and final points of a biarc, starting from $\mathbf{Q}_{k-1}$ and ending in $\mathbf{Q}_{k}$ (see Figure 3 on the right). Moreover, we can easily compute the tangent vectors in $\mathbf{Q}_{k-1}$ and $\mathbf{Q}_{k}$ because we know the centers and the radius of $A_{k-1}$, $A_k$ (a tangent vector of a circle is well-known). The biarc we build with this method regularize the PCC locally, in a neighbourhood of the j.p. $\mathbf{P}_k$. 

It's highly recommended to choose $\delta$ not too high nor too small. In the first case, we can obtain a biarc that doesn't satisfy the tolerance chosen at the beginning of the procedure --- remember that we are trying to build a PCC whose distance from the given points is less than a positive number we decide ---. In the second case we smooth out the curve with such a small radius that the approximating biarc would be very similar to our j.p., i.e. useless. 

A second way to accomplish this task is to use simple arcs: given two arcs (possibly degenerate in segments) with a common j.p., it's always possible to build a simple arc approximating the j.p. within the tolerance and tangent to the other ones. If $R_1$, $R_2$ are the radii of the two consecutive arcs, than we can increase or decrease them to obtain a new circle intersection whose distance from the two arcs is the same. This point will be the center of the tangent arc (see Figure 4).

\section{Ensuring $\mathbf{C^0}$ continuity and reducing the number of arcs}

In the previous sections we discussed about how to find the longest arcs, how to fill the gaps among them and how to smooth the final curve. Now we want to describe a different method for approximating the data set, in order to build a $C^0$ continuous circular curve without any gap. This method tries to find a curve composed by a low number of arcs (not necessarily the minimum). This problem comes from a technical motivation: we want the CNC machine to perform a reduced number of movements. Clearly, minimizing the number of arcs means maximizing the number of approximated points for some arc. In this method we will use again the matrix $M$. The idea is quite simple but efficient. To fix ideas, suppose that $M(1) = 5$. Then, we can build a least square arc starting from $\mathbf{P}_1$, ending in $\mathbf{P}_5$ approximating $\mathbf{P}_2,\mathbf{P}_3,\mathbf{P}_4$. Suppose also that $M(3) = 3$, $M(4) = 7$, $M(5) = 4$. Let $A_1$ be the arc starting from $\mathbf{P}_1$, ending in $\mathbf{P}_5$ and let $A_2$ be the arc starting from $\mathbf{P}_5$, ending in $\mathbf{P}_8$. Then $A_1$, $A_2$ will be $C^0$ continuous in $\mathbf{P}_5$ and will approximate eight points in total. Now let $A_3$ be the arc starting from $\mathbf{P}_1$, ending in $\mathbf{P}_4$ and let $A_4$ be the arc starting from $\mathbf{P}_4$, ending in $\mathbf{P}_{10}$. Then $A_3$, $A_4$ will be $C^0$ continuous in $\mathbf{P}_4$ and will approximate ten points in total. So $A_3$, $A_4$ will be better than $A_1$, $A_2$, because they can approximate a higher number of points. 
Starting from the example we can develop the general case. Let $S = M(i)$, $i = 1,\dots,N-2$, where $N$ is the number of points to be approximated. Then we can find a least square arc starting from $\mathbf{P}_i$, ending in $\mathbf{P}_{i+S-1}$, and a least square arc starting from $\mathbf{P}_{i+S-1}$, ending in $\mathbf{P}_{i+S+M(i+S-1)-2}$. These two arcs will approximate $ S + M(i + S - 1) - 1$ points in total. Now let $j$ be an index such that $i+2 \leq j \leq i + M(i) - 2$. We can build with the same technique two consecutive continuous arcs, the first one from $\mathbf{P}_i$ to $\mathbf{P}_j$ and the second one from $\mathbf{P}_j$ to $\mathbf{P}_{j+M(j)-1}$, approximating $j - i + 1 + M(j) - 1 = j - i + M(j)$ points in total. So, if $j - i + M(j) > S + M(i + S - 1) - 1$, these two arcs will be better than the first ones.
Also here we will work with a new matrix $G$: at the beginning we set $G = M$; during the next steps, we will set $G(j) = 0$ if $G(j)$ is not a starting point for some arcs (as we have already done in the previous sections).

\hrulefill

\begin{enumerate}

\item Set $i = 1$, $G = M$;
\item \textbf{While} $i \leq N-2$ \textbf{do}
\addtolength{\itemindent}{1cm}
\item Set $S = M(i)$, $k = i$, $j = i+2$;
\item Set $Z = S + M(i+S-1)-1$;
\item \textbf{While} $j \leq i + M(i) - 2$ \textbf{do}
\addtolength{\itemindent}{1cm}
\item \textbf{If} $j-i+M(j) > Z$

\hspace{54pt} $k = j$; $Z = j - i + M(j)$;
\item $j = j + 1$;
\addtolength{\itemindent}{-1cm}
\item \textbf{For} $j = i+1,\dots,k-1$

\hspace{26pt} $G(j) = 0$;
\item \textbf{For} $j = k+1,\dots,k+W(k)-2$

\hspace{26pt} $G(j) = 0$;
\item \textbf{If} $i \neq k$

\hspace{26pt} $G(i) = k - i + 1$;
\item $i = k + M(k) - 1$;

\end{enumerate}
\hrulefill

The matrix $G$ determined in the algorithm is used then in the real approximation, and we can apply again the algorithm studied previously. As one can immediately notice, this method is better than the longest-arc-search algorithm if we take in consideration the regularity. Indeed, the final curve will be $C^0$ continuous. The worst case occurs when the final circular spline doesn't end precisely in the last point $\mathbf{P}_N$ but in $\mathbf{P}_{N-1}$. In this situation we have to build a segment from $\mathbf{P}_{N-1}$ to $\mathbf{P}_N$ in order to achieve $C^0$ continuity.

\section{Corner detection}

Up to now we have supposed that all the sections we want to approximate don't contain any corner points. In this paragraph we are going to explain how to detect corner points in the data set. This is not an easy problem, even when the points to be approximated are not affected by errors.
We have to distinguish two different cases: the first one occurs when the data set actually contains corner points, i.e. there are some points that are also corners; the second one occurs when the corners are not in the data set, i.e. all the points are not corners. It's easy to understand that the second case is more difficult and challenging than the first one. This case is caused by detection errors: several laser detectors can't give the exact position of corners. In presence of a corner, they create a sequence of points really close to it, but each of those points can't be considered as a good approximation of the corner. The first case can be solved with one of the algorithms explained in [1], [2], [6], [11], [12]. So, let us focus on the second case. The trivial idea is the following: suppose our data set is composed by $N$ distinct points $\mathbf{P}_i$; for each point we define two unit vectors
\[
\begin{split}
& \mathbf{T}_{1,i} = \frac{\mathbf{P}_{i-1} - \mathbf{P}_i}{\|\mathbf{P}_{i-1} - \mathbf{P}_i\|},\\ 
& \mathbf{T}_{2,i} = \frac{\mathbf{P}_{i+1} - \mathbf{P}_i}{\|\mathbf{P}_{i+1} - \mathbf{P}_i\|},
\end{split}
\]
and we compute the cosine of the angle $\alpha_i$ between the two vectors, using the formula $\cos{\alpha_i} = \mathbf{T}_{1,i} \cdot {\mathbf{T}_{2,i}}$. Then, we say that $\mathbf{P}_i$ is a corner if $\cos{\alpha_i} \geq -\varepsilon$, where $\varepsilon$ is a positive value in $[0,1)$. With few examples one can easily see that this condition is not sufficient to define corners, especially when the data set falls in the second case (when points are not corners). So we need another condition. The second requirement comes from a geometrical interpretation of curvature. We know that the curvature of a smooth curve at each point is the reciprocal of the radius of its osculating circle. If the curvature at the point $x$ is a small number, then the curve will be similar to a line in $x$. 

\begin{figure*}[htbp]
\centering
\includegraphics[scale=0.4]{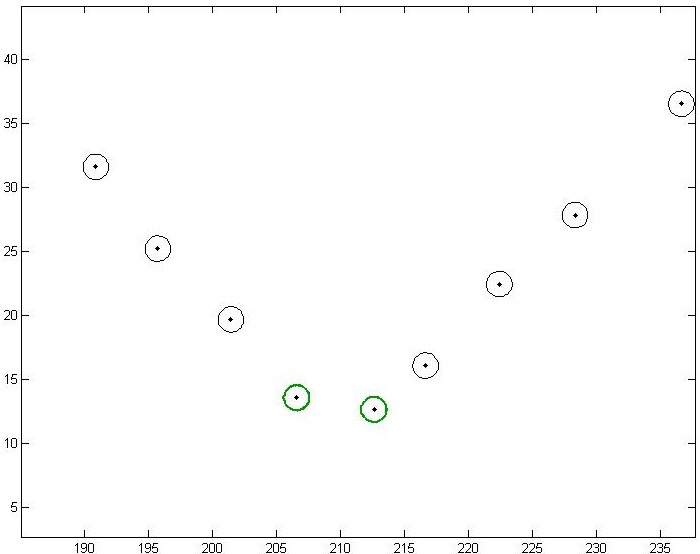}%
\qquad
\includegraphics[scale=0.4]{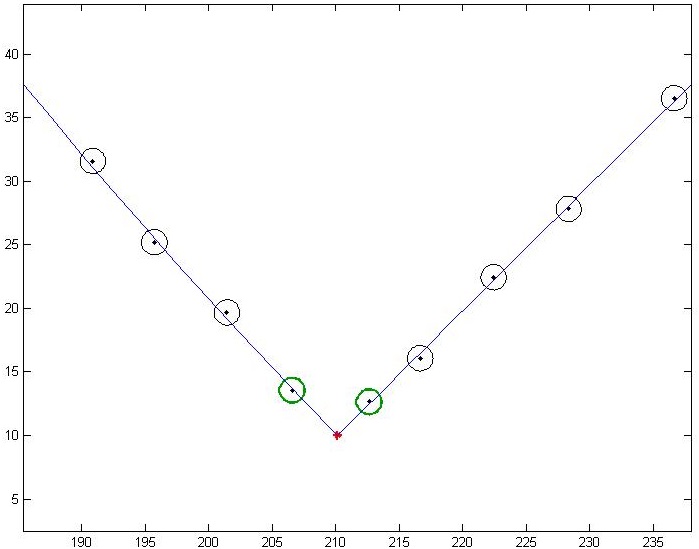}%
\caption{Example of undefined corner (on the left) and approximation of it (red point on the right). Here a tolerance of $\delta = 1\ mm$ was used. The green circles represent the anchor points detected by our algorithm using $R_{max} = 20\ mm$ and $\varepsilon = 0.98$. The blue lines are least square contrained arcs.}
\label{fig: esempi1}
\end{figure*}

Instead, if the curvature at $x$ is high, then the curve will be similar to a small circle in $x$, i.e. the curve will quickly change \textit{direction}. We will use the same principle to determine corner points. 

Assume now that all the points of the data set have been captured by laser detector with constant step. Then, the distance from two consecutive points will be approximately the same. Let us choose two values $\varepsilon \in [0,1)$, $R_{max} > 0$ and consider three consecutive points $\mathbf{P}_{i-1}$, $\mathbf{P}_i$, $\mathbf{P}_{i+1}$ satisfying  $\cos{\alpha_i} \geq -\varepsilon$. Let $R$ be the radius of the circle passing through those three points. We say that $\mathbf{P}_i$ is an \textit{anchor} point if $R < R_{max}$. 

Notice this: being an anchor point doesn't imply being a corner point. If $\mathbf{P}_i$ is an anchor point, then it's quite probable that the actual corner point is in a neighbourhood of  $\mathbf{P}_i$. In order to localize the corner we need two anchor points. Then, we compute $\alpha_{i+1}$, $\alpha_{i-1}$ with the same technique described above. If $\cos{\alpha_{i+1}} \leq \cos{\alpha_{i-1}}$, the two anchor points will be $\mathbf{P}_i$ and $\mathbf{P}_{i+1}$. Otherwise, they will be  $\mathbf{P}_i$ and $\mathbf{P}_{i-1}$. Here one can understand the assumption that the distance from any two consecutive points is approximately the same: if the points are detected with adaptive steps techniques, it's hard to find a right value of $R_{max}$. 

After finding all the anchor points (two at a time), we use them to find corner points. Clearly, the corners need to be placed between these two points. The idea is to use another time the least square arcs. Suppose that $\mathbf{P}_i$ and $\mathbf{P}_{i+1}$ are two consecutive anchor points. Let us consider the sets
\[ 
\begin{split}
&\Lambda_{i,k} = \{ \mathbf{P}_{i-j} : j = 0,\dots,k\},\\
&\Lambda_{i+1,h} = \{ \mathbf{P}_{i+j} : j = 1,\dots,h\},
\end{split}\]
where $k,\ h \in \mathds{N}$ and $k \geq 2,\ h \geq 3$. We want to increase $k,\ h$ in order to build two arcs $A_{i,k}$, $A_{i+1,h}$ approximating respectively $\Lambda_{i,k}$, $\Lambda_{i+1,h}$ within a given positive tolerance $\delta$. Notice this: the arcs we are building are not constrained arcs, i.e. they don't pass through any two given points. This time the best solution comes from Taubin's algorithm. For $k = 2,\ h = 3$, the solution is trivial because the circle passing through three distinct points always exists (degenerate in a line, eventually). So we try to increase $k$ and $h$ (up to a maximum limit $M$) and we build Taubin's circles. If the distance between $\Lambda_{i,k}$ and $A_{i,k}$ is less than $\delta$, $k$ need to be increased.

Otherwise, we decrease $k$ and, in this case, we build the three point circle. We do the same control for $A_{i+1,h}$. At the end of the process, we will find two right values for $k$ and $h$ and we will build the two corresponding arcs $A_{i,k}$, $A_{i+1,h}$. These two arcs must intersect each other in two points. We select the nearest point to $\mathbf{P}_i$. This point will be the approximation of the corner point. 

Pay attention to the fact that if the circles don't intersect each other, then there is something wrong in the data set or the algorithm detected two false anchor points. Another consideration: the corner points will be located approximately between the two anchor points because we used circles to find them: Taubin's method allows to create the best circle approximating a sequence of points, so it's reasonable to assume that $A_i$, $A_{i+1}$ are locally good approximations of the object contour. 

\hrulefill
\begin{enumerate}
\item Find two consecutive anchor points $\mathbf{P}_i$, $\mathbf{P}_{i+1}$;
\item Set $k = 3$, $h = 4$, $flag_k = 1$, $flag_h = 1$ , choose an integer $M \geq 4$ and a tolerance $\delta > 0$;
\item \textbf{While} $k,\ h \leq M$ \textbf{do}
\addtolength{\itemindent}{1cm}
\item Find Taubin's circles $A_{i,k}$, $A_{i+1,h}$ 

\hspace{25pt} approximating $\Lambda_{i,k}$, $\Lambda_{i+1,h}$;
\item \textbf{If} $flag_k = 1$
\addtolength{\itemindent}{1cm}
\item \textbf{If} $\text{dist}(A_{i,k},\Lambda_{i,k}) < \delta$

\hspace{55pt} $k = k + 1$;

\hspace{55pt} \textbf{else} $flag_k = 0$;
\addtolength{\itemindent}{-1cm}
\item \textbf{If} $flag_h = 1$
\addtolength{\itemindent}{1cm}
\item \textbf{If} $\text{dist}(A_{i+1,h},\Lambda_{i+1,h}) < \delta$ 

\hspace{55pt} $h = h + 1$;

\hspace{55pt} \textbf{else} $flag_h = 0$;
\addtolength{\itemindent}{-1cm}
\item \textbf{If} $flag_k = 0$ \textbf{and} $flag_h = 0$

\hspace{26pt} \textbf{return} $k-1$, $h-1$; \textbf{break};

\addtolength{\itemindent}{-1cm}
\item \textbf{Return} $k-1$, $h-1$.
\end{enumerate} 
\hrulefill

\section{Testing the methods}

In this last section we want to give some examples of application of the two methods. First, we want to find a good approximation for the contour of the object represented in Figure 6 on the next page. We give a first rough approximation with 134 points captured by a laser detector (see Figure 6 again). 

\begin{figure*}[htbp]
\centering
\includegraphics[scale=0.19]{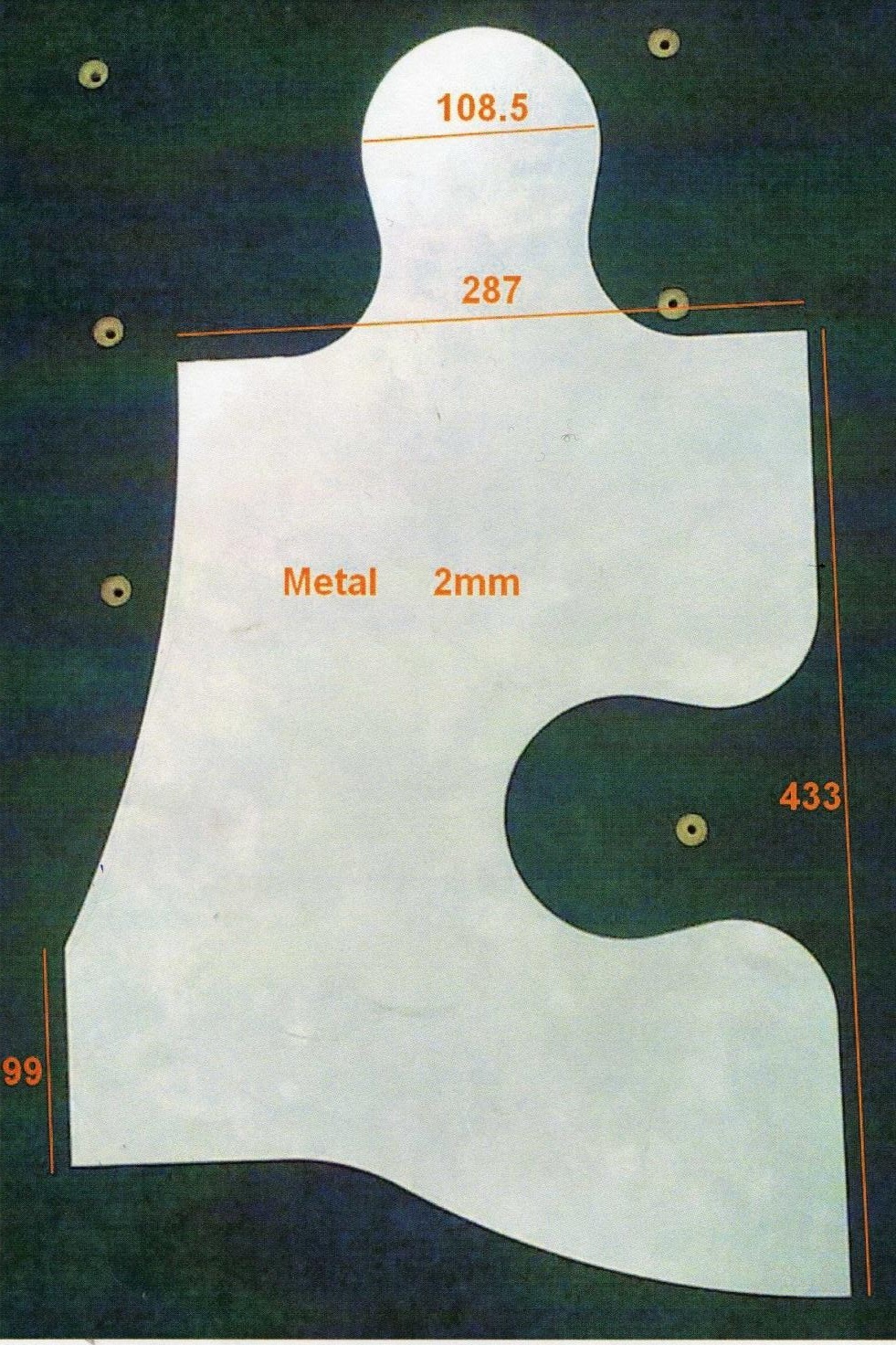}%
\qquad
\includegraphics[scale=0.52]{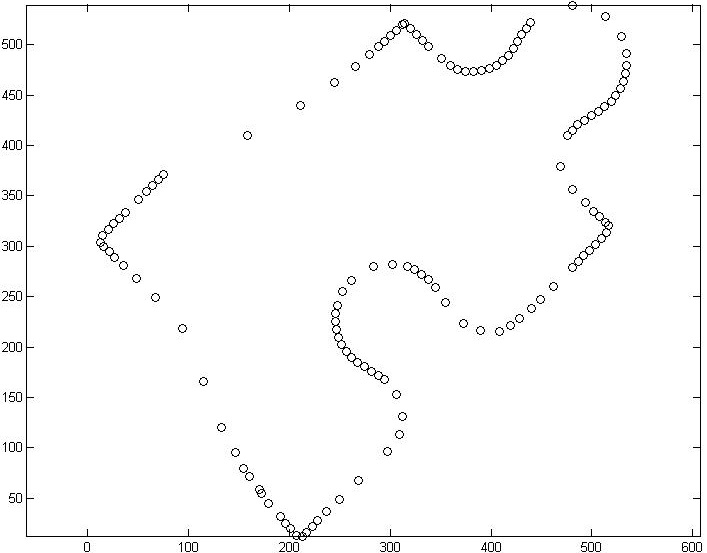}%
\caption{Metal 2D object to be approximated (on the left) and points captured by laser detector (on the right).}
\label{fig: esempi2}
\end{figure*}

\begin{figure*}[htbp]
\centering
\includegraphics[scale=0.40]{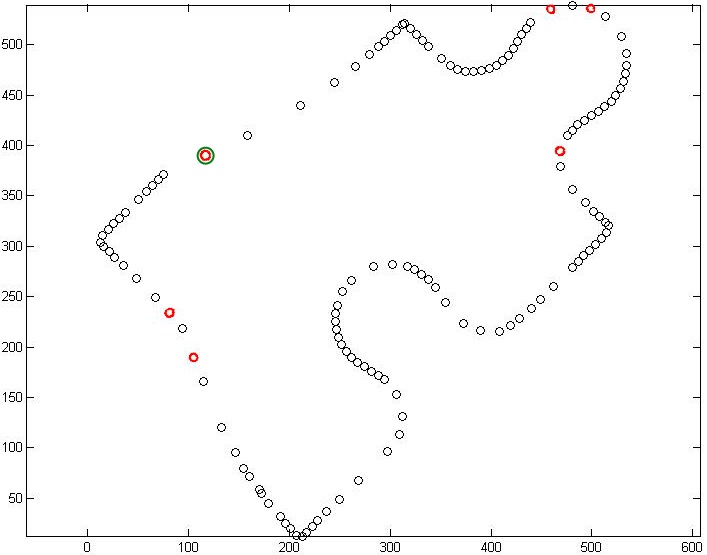}%
\qquad
\includegraphics[scale=0.40]{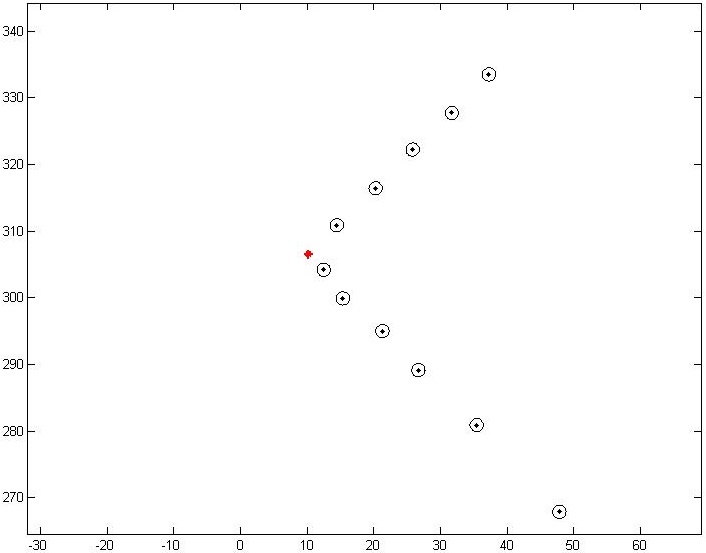}%
\caption{Improvement of the data set.}
\label{fig: esempi3}
\end{figure*}

\begin{figure*}[htbp]
\centering
\includegraphics[scale=0.55]{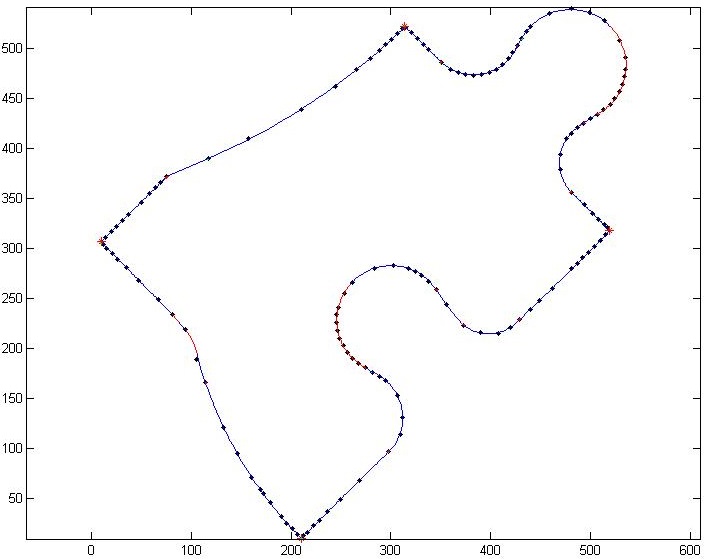}%
\caption{Approximation of the data set using longest arcs approach.}
\label{fig: esempi4}
\end{figure*}

As one can see, this is a bad situation: the points are not uniform, i.e. they were determined with an adaptive method. We select this example because we want to show that the two methods also work in difficult cases. First, we try to improve the data set: since there are wide gaps between some points, we add six points (following a procedure that uses arc intersections) and we obtain the data set represented on Figure 7. Among the new points, the most important one is the green one, because the gap in that section of the contour is too wide to hope for a correct approximation. Secondly, we apply the detection algorithms and we find all five corner points. We use $R_{max} = 20 mm$, $\varepsilon = 0.9$ (for detection of anchor points). 

\begin{table}[ht]
\caption{Results with the first method.}
\begin{center}
\begin{tabular}{cccc}
\toprule
Tolerance & Arcs & Segments & Bad Points \\
\midrule
1.5 mm & 14 & 1 & 2\\
1 mm & 17 & 3 & 6\\
0.5  mm & 22 & 3 & 6\\
\bottomrule
\label{tab: one}
\end{tabular}
\end{center}
\end{table}

\begin{table}[ht]
\caption{Results with the second method.}
\begin{center}
\begin{tabular}{cccc}
\toprule
Tolerance & Arcs & Segments & Bad Points \\
\midrule
1.5 mm & 18 & 0 & 3\\
1 mm & 21 & 1 & 2\\
0.5  mm & 28 & 0 & 4\\
\bottomrule
\label{tab: two}
\end{tabular}
\end{center}
\end{table}

Notice this: adding the green points is crucial. Indeed, in this case the research of all the corner points fails, because the initial data set is ill-posed in a neighbourhood of one of the corners, i.e. there are not enough points to define a good approximation. Using the two methods, we find the results in Table\footnote{Here the tolerance is expressed in millimetres. Remember that 1 mm $\approx$ 0.03937 Inches} 1 and 2. As we expected, the number of arcs in the first method is less than the number of arcs in the second one. But using the second approach we create a continuous curve, so the number of segments is reduced. 

\begin{figure*}[htbp]
\centering
\includegraphics[scale=0.40]{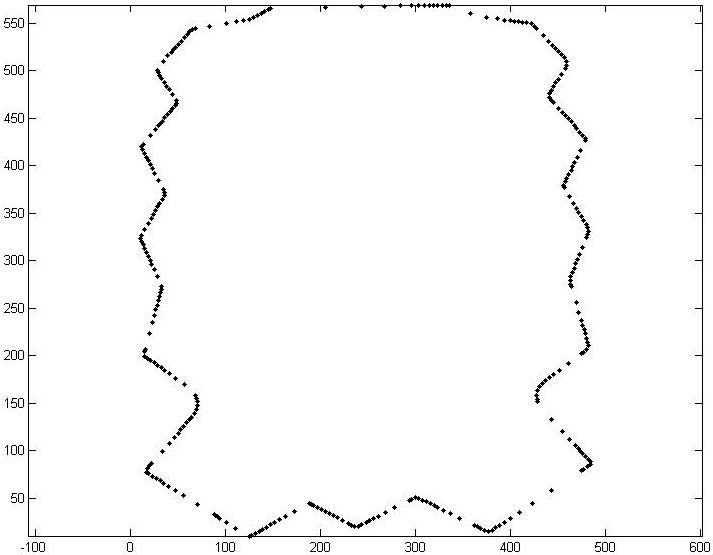}%
\qquad
\includegraphics[scale=0.40]{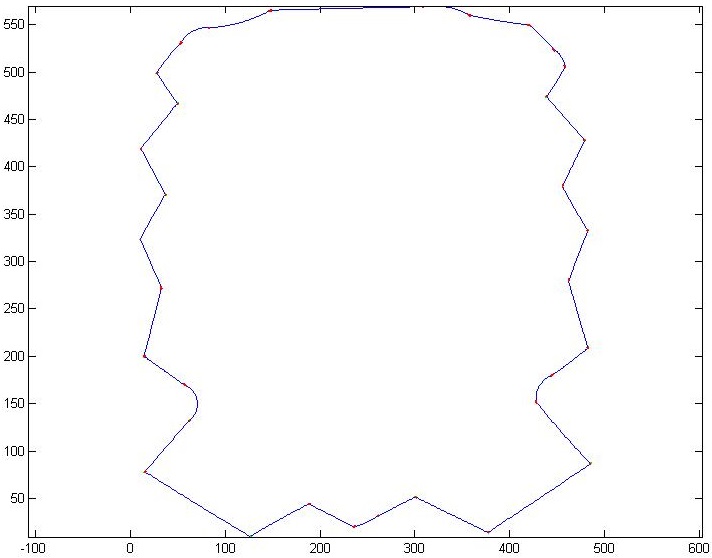}%
\caption{Contour of a wooden object (left) and approximation of it (right).}
\label{fig: scan_wood}
\end{figure*}

\begin{figure*}[htbp]
\centering
\includegraphics[scale=0.40]{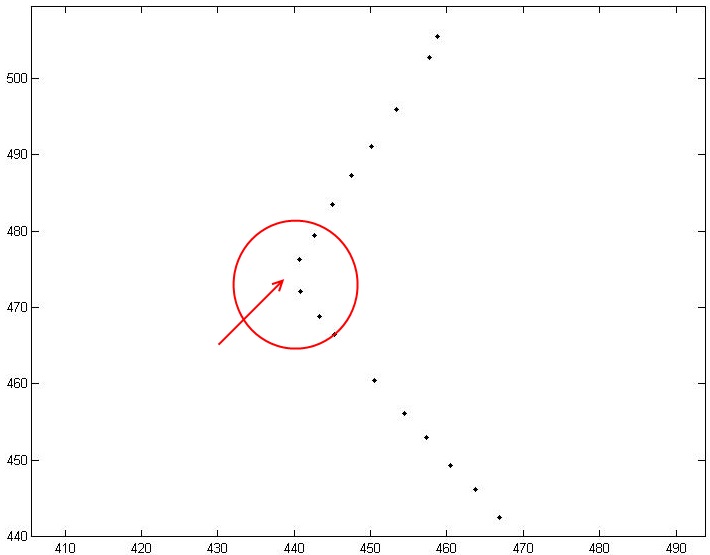}%
\qquad
\includegraphics[scale=0.40]{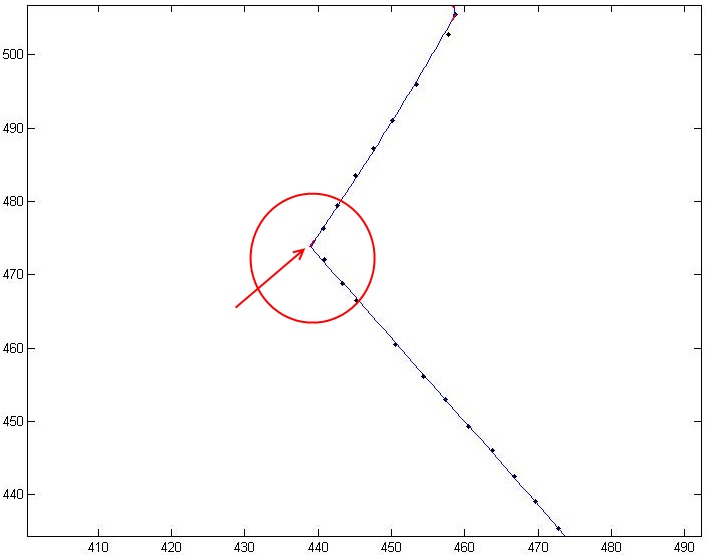}%
\caption{The true corners are not included in the data set.}
\label{fig: scan_wood_corners}
\end{figure*}

Notice also that this method works better if the points are detected precisely. In our example, there are some detection errors, so it is better to use large error tolerance. If we try to use a tolerance smaller than $0.1$, the final result is not reliable, because in the best case we would obtain a curve composed by arcs approximating a low number of points (three or four), so we would also have a high number of bad points.

Secondly, we want to approximate the contour of Figure 9. This time we don't need to build new points, because the sections between two consecutive corners are quite short and linear. However, we need to find a good approximation for corner points: as we can see in Figure 10 there is the same problem we faced in the other example, so the corner points are not included in the data set. Using $R_{max} = 20mm$ and $\varepsilon = 0.85$ (for detection of anchor points) we find 18 corner points. After applying the second method, we have the results in Table 3.

\begin{table}[ht]
\caption{Results with the second method on the wooden object with $\varepsilon = 0.85$.}
\begin{center}
\begin{tabular}{cccc}
\toprule
Tolerance & Arcs & Segments & Bad Points \\
\midrule
1.5 mm & 43 & 0 & 9\\
1 mm & 46 & 0 & 10\\
0.5  mm & 51 & 0 & 10\\
\bottomrule
\label{tab: three}
\end{tabular}
\end{center}
\end{table}

\begin{table}[ht]
\caption{Results with the second method on the wooden object with $\varepsilon = 0.9$.}
\begin{center}
\begin{tabular}{cccc}
\toprule
Tolerance & Arcs & Segments & Bad Points \\
\midrule
1.5 mm & 38 & 0 & 4\\
1 mm & 40 & 0 & 7\\
0.5  mm & 43 & 0 & 7\\
\bottomrule
\label{tab: four}
\end{tabular}
\end{center}
\end{table}

We can notice that the number of bad points in Table 3 is quite high. This is justified by the fact that the algorithm detected only 18 corners, so at least three or four bad points should be corner points. Indeed, using $\varepsilon = 0.9$ for detection of anchor points we find 21 corners and the results in Table 4.

Finally, we want to test the method in a simple case. The contour to be approximated is represented in Figure 11.

\begin{figure}[htbp]
\centering
\includegraphics[scale=0.40]{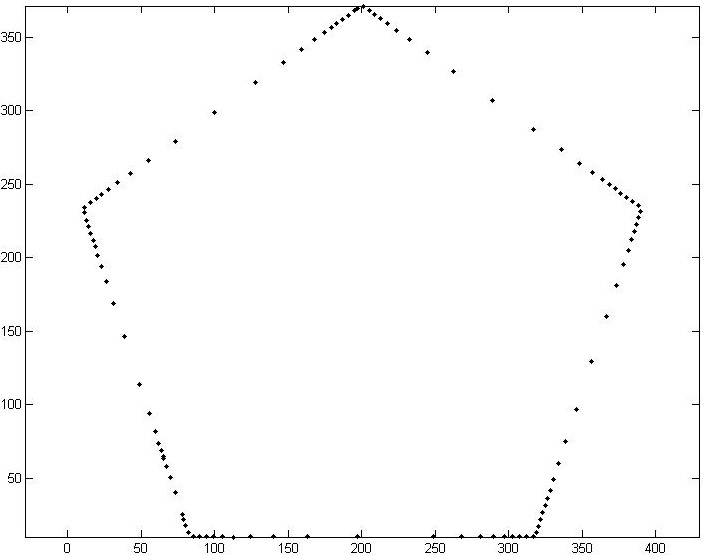}%
\caption{Simple contour.}
\label{fig: pentagono}
\end{figure}

In this situation, the sections between two consecutive corners are simple lines, but there is still the problem of detecting the true corners, as we can notice in Figure 12.
Using again $R_{max} = 20mm$ and $\varepsilon = 0.9$ we find five corners and the results of Table 5.

\begin{table}[ht]
\caption{Results with the second method on the simple contour.}
\begin{center}
\begin{tabular}{cccc}
\toprule
Tolerance & Arcs & Segments & Bad Points \\
\midrule
1.5 mm & 5 & 0 & 0\\
1 mm & 5 & 0 & 0\\
0.5  mm & 5 & 0 & 0\\
\bottomrule
\label{tab: five}
\end{tabular}
\end{center}
\end{table}

\begin{figure*}[htbp]
\centering
\includegraphics[scale=0.40]{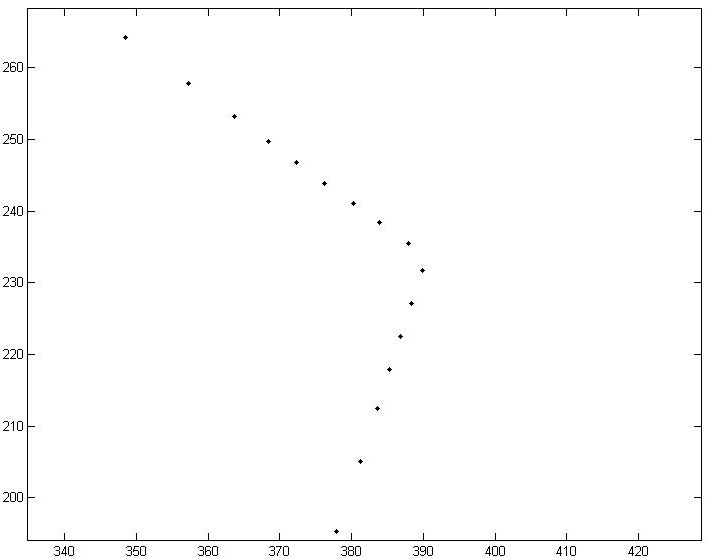}%
\qquad
\includegraphics[scale=0.40]{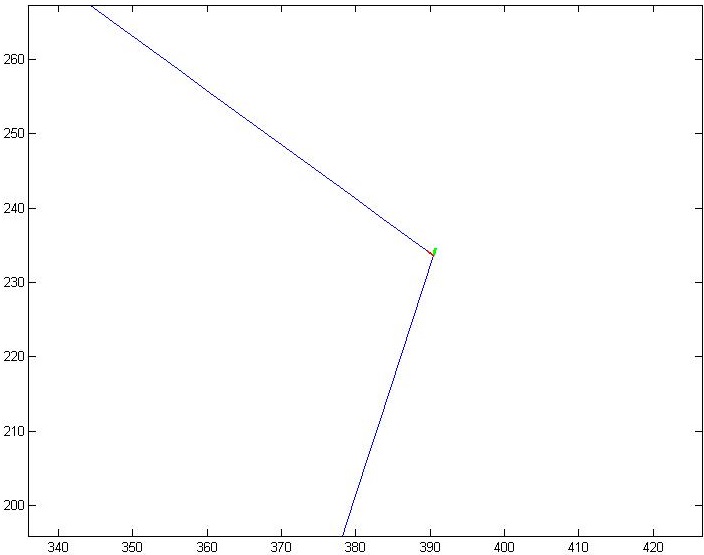}%
\caption{Also in this case the corners are not included in the data set.}
\label{fig: pentagono1}
\end{figure*}

\section*{Conclusions}
Two methods of approximation of a planar data set were presented. Both methods are based on arcs and segments. The first one searches the longest arcs in the data set and eventually fills the gaps among them using biarcs or segments. The second one creates a lower number of arcs but the final curve is more regular and $C^0$ continuous. Moreover, the number of segments is minimized. Finally, a method of corner detection was presented. The corners of the object contour are determined using circle intersections and Taubin's algorithm. 

\section*{Aknowledgments}
The author is very greatful to Professor C. Dagnino and Dr. D. Inaudi for their support and their suggestions.


%
%

\end{document}